# Fractal properties of the random string processes

Dongsheng Wu[1],[*] and Yimin Xiao[1],[†]

*Department of Statistics and Probability, Michigan State University*

**Abstract:** Let $\{u_t(x), t \geq 0, x \in \mathbb{R}\}$ be a random string taking values in $\mathbb{R}^d$, specified by the following stochastic partial differential equation [Funaki (1983)]:
$$\frac{\partial u_t(x)}{\partial t} = \frac{\partial^2 u_t(x)}{\partial x^2} + \dot{W},$$
where $\dot{W}(x,t)$ is an $\mathbb{R}^d$-valued space-time white noise.

Mueller and Tribe (2002) have proved necessary and sufficient conditions for the $\mathbb{R}^d$-valued process $\{u_t(x) : t \geq 0, x \in \mathbb{R}\}$ to hit points and to have double points. In this paper, we continue their research by determining the Hausdorff and packing dimensions of the level sets and the sets of double times of the random string process $\{u_t(x) : t \geq 0, x \in \mathbb{R}\}$. We also consider the Hausdorff and packing dimensions of the range and graph of the string.

## 1. Introduction and preliminaries

Consider the following model of a random string introduced by Funaki [5]:

(1) $$\frac{\partial u_t(x)}{\partial t} = \frac{\partial^2 u_t(x)}{\partial x^2} + \dot{W},$$

where $\dot{W}(x,t)$ is a space-time white noise in $\mathbb{R}^d$, which is assumed to be adapted with respect to a filtered probability space $(\Omega, \mathcal{F}, \mathcal{F}_t, \mathbb{P})$, where $\mathcal{F}$ is complete and the filtration $\{\mathcal{F}_t, t \geq 0\}$ is right continuous. The components $\dot{W}_1(x,t), \ldots, \dot{W}_d(x,t)$ of $\dot{W}(x,t)$ are independent space-time white noises, which are generalized Gaussian processes with covariance given by

$$\mathbb{E}\big[\dot{W}_j(x,t)\dot{W}_j(y,s)\big] = \delta(x-y)\delta(t-s), \quad (j = 1,\ldots,d).$$

That is, for every $1 \leq j \leq d$, $W_j(f)$ is a random field indexed by functions $f \in L^2([0,\infty) \times \mathbb{R})$ and, for all $f, g \in L^2([0,\infty) \times \mathbb{R})$, we have

$$\mathbb{E}\big[W_j(f)W_j(g)\big] = \int_0^\infty \int_\mathbb{R} f(t,x)g(t,x)\,dxdt.$$

Hence $W_j(f)$ can be represented as

$$W_j(f) = \int_0^\infty \int_\mathbb{R} f(t,x)\,W_j(dx\,dt).$$

---

[1]Department of Statistics and Probability, Michigan State University, East Lansing, MI 48824, USA, e-mail: wudongsh@msu.edu; xiao@stt.msu.edu; url: www.msu.edu/~wudongsh; www.stt.msu.edu/~xiaoyimi

[*]Research partially supported by NSF grant DMS-0417869.
[†]Research partially supported by NSF grant DMS-0404729.

Note that $W(f)$ is $\mathcal{F}_t$-measurable whenever $f$ is supported on $[0, t] \times \mathbb{R}$.

Recall from Mueller and Tribe [9] that a solution of (1) is defined as an $\mathcal{F}_t$-adapted, continuous random field $\{u_t(x) : t \geq 0, x \in \mathbb{R}\}$ with values in $\mathbb{R}^d$ satisfying the following properties:

(i) $u_0(\cdot) \in \mathcal{E}_{\exp}$ almost surely and is adapted to $\mathcal{F}_0$, where $\mathcal{E}_{\exp} = \cup_{\lambda > 0} \mathcal{E}_\lambda$ and
$$\mathcal{E}_\lambda = \left\{ f \in C(\mathbb{R}, \mathbb{R}^d) : |f(x)| e^{-\lambda |x|} \to 0 \text{ as } |x| \to \infty \right\};$$

(ii) For every $t > 0$, there exists $\lambda > 0$ such that $u_s(\cdot) \in \mathcal{E}_\lambda$ for all $s \leq t$, almost surely;

(iii) For every $t > 0$ and $x \in \mathbb{R}$, the following Green's function representation holds
$$(2) \qquad u_t(x) = \int_\mathbb{R} G_t(x - y) u_0(y) dy + \int_0^t G_{t-r}(x - y) W(dy\, dr),$$

where $G_t(x) = \frac{1}{\sqrt{4\pi t}} e^{-\frac{x^2}{4t}}$ is the fundamental solution of the heat equation.

We call each solution $\{u_t(x) : t \geq 0, x \in \mathbb{R}\}$ of (1) a random string process with values in $\mathbb{R}^d$, or simply a random string as in [9]. Note that, whenever the initial conditions $u_0$ are deterministic, or are Gaussian fields independent of $\mathcal{F}_0$, the random string processes are Gaussian. We recall briefly some basic properties about the solutions of (1), and refer to Mueller and Tribe [9] and Funaki [5] for further information on stochastic partial differential equations (SPDEs) related to random motion of strings.

Funaki [5] investigated various properties of the solutions of semi-linear type SPDEs which are more general than (1). In particular, his results (cf. Lemma 3.3 in [5]) imply that every solution $\{u_t(x) : t \geq 0, x \in \mathbb{R}\}$ of (1) is Hölder continuous of any order less than $\frac{1}{2}$ in space and $\frac{1}{4}$ in time. This anisotropic property of the process $\{u_t(x) : t \geq 0, x \in \mathbb{R}\}$ makes it a very interesting object to study. Recently Mueller and Tribe [9] have found necessary and sufficient conditions [in terms of the dimension $d$] for a random string in $\mathbb{R}^d$ to hit points or to have double points of various types. They have also studied the question of recurrence and transience for $\{u_t(x) : t \geq 0, x \in \mathbb{R}\}$. Note that, in general, a random string may not be Gaussian, a powerful step in the proofs of Mueller and Tribe [9] is to reduce the problems about a general random string process to those of the stationary pinned string $U = \{U_t(x), t \geq 0, x \in \mathbb{R}\}$, obtained by taking the initial functions $U_0(\cdot)$ in (2) to be defined by

$$(3) \qquad U_0(x) = \int_0^\infty \int (G_r(x - z) - G_r(z)) \widetilde{W}(dz dr),$$

where $\widetilde{W}$ is a space-time white noise independent of the white noise $\dot{W}$. One can verify that $U_0 = \{U_0(x) : x \in \mathbb{R}\}$ is a two-sided $\mathbb{R}^d$ valued Brownian motion satisfying $U_0(0) = 0$ and $\mathbb{E}[(U_0(x) - U_0(y))^2] = |x - y|$. We assume, by extending the probability space if needed, that $U_0$ is $\mathcal{F}_0$-measurable. As pointed out by Mueller and Tribe [9], the solution to (1) driven by the noise $W(x, s)$ is then given by

$$(4) \qquad \begin{aligned} U_t(x) &= \int G_t(x - z) U_0(z) dz + \int_0^t \int G_r(x - z) W(dz dr) \\ &= \int_0^\infty (G_{t+r}(x - z) - G_r(z)) \widetilde{W}(dz dr) + \int_0^t \int G_r(x - z) W(dz dr). \end{aligned}$$

A continuous version of the above solution is called *a stationary pinned string*.

The components $\{U_t^j(x) : t \geq 0, x \in \mathbb{R}\}$ for $j = 1, \ldots, d$ are independent and identically distributed Gaussian processes. In the following we list some basic properties of the processes $\{U_t^j(x) : t \geq 0, x \in \mathbb{R}\}$, which will be needed for proving the results in this paper. Lemma 1.1 below is Proposition 1 of Mueller and Tribe [9].

**Lemma 1.1.** *The components $\{U_t^j(x) : t \geq 0, x \in \mathbb{R}\}$ ($j = 1, \ldots, d$) of the stationary pinned string are mean-zero Gaussian random fields with stationary increments. They have the following covariance structure: for $x, y \in \mathbb{R}$, $t \geq 0$,*

$$\mathbb{E}\Big[\big(U_t^j(x) - U_t^j(y)\big)^2\Big] = |x - y|, \tag{5}$$

*and for all $x, y \in \mathbb{R}$ and $0 \leq s < t$,*

$$\mathbb{E}\Big[\big(U_t^j(x) - U_s^j(y)\big)^2\Big] = (t-s)^{1/2} F\big(|x-y|(t-s)^{-1/2}\big), \tag{6}$$

*where*

$$F(a) = (2\pi)^{-1/2} + \frac{1}{2}\int_{\mathbb{R}}\int_{\mathbb{R}} G_1(a-z)G_1(a-z')\big(|z| + |z'| - |z - z'|\big)\,dz\,dz'.$$

$F(x)$ *is a smooth function, bounded below by $(2\pi)^{-1/2}$, and $F(x)/|x| \to 1$ as $|x| \to \infty$. Furthermore there exists a positive constant $c_{1,1}$ such that for all $s, t \in [0, \infty)$ and all $x, y \in \mathbb{R}$,*

$$c_{1,1}\big(|x-y| + |t-s|^{1/2}\big) \leq \mathbb{E}\Big[\big(U_t^j(x) - U_s^j(y)\big)^2\Big] \leq 2\big(|x-y| + |t-s|^{1/2}\big). \tag{7}$$

It follows from (6) that the stationary pinned string has the following scaling property [or operator-self-similarity]: For any constant $c > 0$,

$$\{c^{-1}U_{c^4 t}(c^2 x) : t \geq 0, x \in \mathbb{R}\} \stackrel{d}{=} \{U_t(x) : t \geq 0, x \in \mathbb{R}\}, \tag{8}$$

where $\stackrel{d}{=}$ means equality in finite dimensional distributions; see Corollary 1 in [9].

We will also need more precise information about the asymptotic property of the function $F(x)$. By a change of variables we can write it as

$$F(x) = -(2\pi)^{-1/2} + \frac{1}{2}\int_{\mathbb{R}}\int_{\mathbb{R}} G_1(z)G_1(z')\big(|z - x| + |z' - x|\big)\,dz\,dz'. \tag{9}$$

Denote the above double integral by $H(x)$. Then it can be written as

$$H(x) = \int_{\mathbb{R}} G_1(z)\,|z - x|\,dz. \tag{10}$$

The following lemma shows that the behavior of $H(x)$ is similar to that of $F(x)$, and the second part describes how fast $H(x)/|x| \to 1$ as $x \to \infty$.

**Lemma 1.2.** *There exist positive constants $c_{1,2}$ and $c_{1,3}$ such that*

$$c_{1,2}\big(|x-y|+|t-s|^{1/2}\big) \leq |t-s|^{1/2} H\big(|x-y||t-s|^{-1/2}\big) \leq c_{1,3}\big(|x-y|+|t-s|^{1/2}\big). \tag{11}$$

*Moreover, we have the limit:*

$$\lim_{x \to \infty} |H(x) - x| = 0. \tag{12}$$

*Proof.* The inequality (11) follows from the proof of (7) in [9], p. 9. Hence we only need to prove (12).

By (10), we see that for $x > 0$,

$$
\begin{aligned}
H(x) - x &= \int_{\mathbb{R}} G_1(z)\bigl(|z-x| - x\bigr) dz \\
&= \int_x^\infty (z - 2x)\, G_1(z)\, dz - \int_{-\infty}^x z\, G_1(z)\, dz \\
&= 2\int_x^\infty (z - x)\, G_1(z)\, dz.
\end{aligned}
\tag{13}
$$

Since the last integral tends to 0 as $x \to \infty$, (12) follows. □

The following lemmas indicate that, for every $j \in \{1, 2, \ldots, d\}$, the Gaussian process $\{U_t^j(x), t \geq 0, x \in \mathbb{R}\}$ satisfies some preliminary forms of sectorial local nondeterminism; see [13] for more information on the latter. Lemma 1.3 is implied by the proof of Lemma 3 in [9], p. 15, and Lemma 1.4 follows from the proof of Lemma 4 in [9], p. 21.

**Lemma 1.3.** *For any given $\varepsilon \in (0,1)$, there exists a positive constant $c_{1,4}$, which depends on $\varepsilon$ only, such that*

$$
\operatorname{Var}\left(U_t^j(x)\,\Big|\,U_s^j(y)\right) \geq c_{1,4}\left(|x-y| + |t-s|^{1/2}\right)
\tag{14}
$$

*for all $(t,x), (s,y) \in [\varepsilon, \varepsilon^{-1}] \times [-\varepsilon^{-1}, \varepsilon^{-1}]$.*

**Lemma 1.4.** *For any given constants $\varepsilon \in (0,1)$ and $L > 0$, there exists a constant $c_{1,5} > 0$ such that*

$$
\begin{aligned}
&\operatorname{Var}\left(U_{t_2}^j(x_2) - U_{t_1}^j(x_1)\,\Big|\,U_{s_2}^j(y_2) - U_{s_1}^j(y_1)\right) \\
&\geq c_{1,5}\left(|x_1 - y_1| + |x_2 - y_2| + |t_1 - s_1|^{1/2} + |t_2 - s_2|^{1/2}\right)
\end{aligned}
\tag{15}
$$

*for all $(t_k, x_k), (s_k, y_k) \in [\varepsilon, \varepsilon^{-1}] \times [-\varepsilon^{-1}, \varepsilon^{-1}]$, where $k \in \{1,2\}$, such that $|t_2 - t_1| \geq L$ and $|s_2 - s_1| \geq L$.*

Note that in Lemma 1.4, the pairs $t_1$ and $t_2$, $s_1$ and $s_2$, are well separated. The following lemma is concerned with the case when $t_1 = t_2$ and $s_1 = s_2$.

**Lemma 1.5.** *Let $\varepsilon \in (0,1)$ and $L > 0$ be given constants. Then there exist positive constants $h_0 \in (0, \frac{L}{2})$ and $c_{1,6}$ such that*

$$
(16)\quad \operatorname{Var}\left(U_t^j(x_2) - U_t^j(x_1)\,\Big|\,U_s^j(y_2) - U_s^j(y_1)\right) \geq c_{1,6}\left(|s-t|^{1/2} + |x_1 - y_1| + |x_2 - y_2|\right)
$$

*for all $s, t \in [\varepsilon, \varepsilon^{-1}]$ with $|s-t| \leq h_0$ and all $x_k, y_k \in [-\varepsilon^{-1}, \varepsilon^{-1}]$, where $k \in \{1,2\}$, such that $|x_2 - x_1| \geq L$, $|y_2 - y_1| \geq L$ and $|x_k - y_k| \leq \frac{L}{2}$ for $k = 1, 2$.*

**Remark 1.6.** Note that, in the above, it is essential to only consider those $s, t \in [\varepsilon, \varepsilon^{-1}]$ such that $|s-t|$ is small. Otherwise (16) does not hold as indicated by (5). In this sense, Lemma 1.5 is more restrictive than Lemma 1.4. But it is sufficient for the proof of Theorem 4.2.

*Proof.* Using the notation similar to that in [9], we let $(X, Y) = \big(U_t^j(x_2) - U_t^j(x_1), U_s^j(y_2) - U_s^j(y_1)\big)$ and write $\sigma_X^2 = \mathbb{E}(X^2)$, $\sigma_Y^2 = \mathbb{E}(Y^2)$ and $\rho_{X,Y}^2 = \mathbb{E}\big[(X - Y)^2\big]$. Recall that, for the Gaussian vector $(X, Y)$, we have

$$\text{Var}(X|Y) = \frac{\big(\rho_{X,Y}^2 - (\sigma_X - \sigma_Y)^2\big)\big((\sigma_X + \sigma_Y)^2 - \rho_{X,Y}^2\big)}{4\sigma_Y^2}. \tag{17}$$

Lemma 1.1 and the separation condition on $x_k$ and $y_k$ imply that both $\sigma_X^2$ and $\sigma_Y^2$ are bounded from above and below by positive constants. Similar to the proofs of Lemmas 3 and 4 in [9], we only need to derive a suitable lower bound for $\rho_{X,Y}^2$. By using the identity

$$(a - b + c - d)^2 = (a - b)^2 + (c - d)^2 + (a - d)^2 + (b - c)^2 - (a - c)^2 - (b - d)^2$$

and (5) we have

$$\begin{aligned}\rho_{X,Y}^2 &= |t - s|^{1/2} F\big(|x_2 - y_2||t - s|^{-1/2}\big) + |t - s|^{1/2} F\big(|y_1 - x_1||t - s|^{-1/2}\big) \\ &\quad + |x_2 - x_1| - |t - s|^{1/2} F\big(|x_2 - y_1||t - s|^{-1/2}\big) \\ &\quad + |y_1 - y_2| - |t - s|^{1/2} F\big(|x_1 - y_2||t - s|^{-1/2}\big).\end{aligned} \tag{18}$$

By (9), we can rewrite the above equation as

$$\begin{aligned}\rho_{X,Y}^2 &= |t - s|^{1/2} H\big(|x_2 - y_2||t - s|^{-1/2}\big) \\ &\quad + |t - s|^{1/2} H\big(|y_1 - x_1||t - s|^{-1/2}\big) \\ &\quad + |x_2 - x_1| - |t - s|^{1/2} H\big(|x_2 - y_1||t - s|^{-1/2}\big) \\ &\quad + |y_1 - y_2| - |t - s|^{1/2} H\big(|x_1 - y_2||t - s|^{-1/2}\big).\end{aligned} \tag{19}$$

Denote the algebraic sum of the last four terms in (19) by $S$ and we need to derive a lower bound for it. Note that, under the conditions of our lemma, $|x_2 - y_1| \geq \frac{L}{2}$ and $|x_1 - y_2| \geq \frac{L}{2}$. Hence Lemma 1.2 implies that, for any $0 < \delta < \frac{c_{1,2}}{2}$, there exists a constant $h_0 \in (0, \frac{L}{2})$ such that

$$|t - s|^{1/2} H\big(|x_2 - y_1||t - s|^{-1/2}\big) \leq |x_2 - y_1| + \frac{\delta}{2}|t - s|^{1/2} \tag{20}$$

whenever $|t - s| \leq h_0$; and the same inequality holds when $|x_2 - y_1|$ is replaced by $|x_1 - y_2|$. It follows that

$$\begin{aligned}S &\geq \big(|x_2 - x_1| - |x_2 - y_1| + |y_1 - y_2| - |x_1 - y_2|\big) - \delta |t - s|^{1/2} \\ &= -\delta |t - s|^{1/2},\end{aligned} \tag{21}$$

because the sum of the four terms in the parentheses equals 0 under the separation condition. Combining (19) and (11) yields

$$\rho_{X,Y}^2 \geq \frac{c_{1,2}}{2}\big(|t - s|^{1/2} + |x_1 - y_1| + |x_2 - y_2|\big) \tag{22}$$

whenever $x_k$, $y_k$ ($k = 1, 2$) satisfy the above conditions.

By (5), we have $(\sigma_X - \sigma_Y)^2 \leq c\big(|y_1 - x_1| + |x_2 - y_2|\big)^2$. It follows from (17) and (22) that (16) holds whenever $|y_1 - x_1| + |x_2 - y_2|$ is sufficiently small. Finally, a continuity argument as in [9], p. 15 removes this last restriction. This finishes the proof of Lemma 1.5. □

The present paper is a continuation of the paper of Mueller and Tribe [9]. Our objective is to study the fractal properties of various random sets generated by the random string processes. In Section 2, we determine the Hausdorff and packing dimensions of the range $u([0,1]^2)$ and the graph $\mathrm{Gr}u([0,1]^2)$. We also consider the Hausdorff dimension of the range $u(E)$, where $E \subseteq [0,\infty) \times \mathbb{R}$ is an arbitrary Borel set. In Section 3, we consider the existence of the local times of the random string process and determine the Hausdorff and packing dimensions of the level set $L_{\mathbf{u}} = \{(t,x) \in (0,\infty) \times \mathbb{R} : u_t(x) = \mathbf{u}\}$, where $\mathbf{u} \in \mathbb{R}^d$. Finally, we conclude our paper by determining the Hausdorff and packing dimensions of the sets of two kinds of double times of the random string in Section 4.

## 2. Dimension results of the range and graph

In this section, we study the Hausdorff and packing dimensions of the range $u([0,1]^2) = \{u_t(x) : (t,x) \in [0,1]^2\} \subset \mathbb{R}^d$ and the graph $\mathrm{Gr}u([0,1]^2) = \{((t,x), u_t(x)) : (t,x) \in [0,1]^2\} \subset \mathbb{R}^{2+d}$. We refer to Falconer [4] for the definitions and properties of Hausdorff dimension $\dim_{\mathrm{H}}(\cdot)$ and packing dimension $\dim_{\mathrm{P}}(\cdot)$.

**Theorem 2.1.** *Let $\{u_t(x) : t \geq 0, x \in \mathbb{R}\}$ be a random string process taking values in $\mathbb{R}^d$. Then with probability 1,*

$$\dim_{\mathrm{H}} u([0,1]^2) = \min\{d; 6\} \tag{23}$$

*and*

$$\dim_{\mathrm{H}} \mathrm{Gr}u([0,1]^2) = \begin{cases} 2 + \frac{3}{4}d & \text{if } 1 \leq d < 4, \\ 3 + \frac{1}{2}d & \text{if } 4 \leq d < 6, \\ 6 & \text{if } 6 \leq d. \end{cases} \tag{24}$$

*Proof.* Corollary 2 of Mueller and Tribe [9] states that the distributions of $\{u_t(x) : t \geq 0, x \in \mathbb{R}\}$ and the stationary pinned string $U = \{U_t(x) : t \geq 0, x \in \mathbb{R}\}$ are mutually absolutely continuous. Hence it is enough for us to prove (23) and (24) for the stationary pinned string $U = \{U_t(x) : t \geq 0, x \in \mathbb{R}\}$. This is similar to the proof of Theorem 4 of Ayache and Xiao [2]. We include a self-contained proof for reader's convenience.

As usual, the proof is divided into proving the upper and lower bounds separately. For the upper bound in (23), we note that clearly $\dim_{\mathrm{H}} U([0,1]^2) \leq d$ a.s., so we only need to prove the following inequality:

$$\dim_{\mathrm{H}} U([0,1]^2) \leq 6 \qquad \text{a.s.} \tag{25}$$

Because of Lemma 1.1, one can use the standard entropy method for estimating the tail probabilities of the supremum of a Gaussian process to establish the modulus of continuity of $U = \{U_t(x) : t \geq 0, x \in \mathbb{R}\}$. See, for example, Kôno [8]. It follows that, for any constants $0 < \gamma_1 < \gamma_1' < 1/4$ and $0 < \gamma_2 < \gamma_2' < 1/2$, there exist a random variable $A > 0$ of finite moments of all orders and an event $\Omega_1$ of probability 1 such that for all $\omega \in \Omega_1$,

$$\sup_{(s,y),(t,x)\in[0,1]^2} \frac{|U_s(y,\omega) - U_t(x,\omega)|}{|s-t|^{\gamma_1'} + |x-y|^{\gamma_2'}} \leq A(\omega). \tag{26}$$

Let $\omega \in \Omega_1$ be fixed and then suppressed. For any integer $n \geq 2$, we divide $[0,1]^2$ into $n^6$ sub-rectangles $\{R_{n,i}\}$ with sides parallel to the axes and side-lengths $n^{-4}$ and $n^{-2}$, respectively. Then $U([0,1]^2)$ can be covered by the sets $U(R_{n,i})$ ($1 \leq i \leq n^6$). By (26), we see that the diameter of the image $U(R_{n,i})$ satisfies

$$\text{diam} U(R_{n,i}) \leq c_{2,1}\, n^{-1+\delta}, \tag{27}$$

where $\delta = \max\{1 - 4\gamma_1', 1 - 2\gamma_2'\}$. We choose $\gamma_1' \in (\gamma_1, 1/4)$ and $\gamma_2' \in (\gamma_2, 1/2)$ such that

$$(1-\delta)\left(\frac{1}{\gamma_1} + \frac{1}{\gamma_2}\right) > 6.$$

Hence, for $\gamma = \frac{1}{\gamma_1} + \frac{1}{\gamma_2}$, it follows from (27) that

$$\sum_{i=1}^{n^6} \left[\text{diam} U(R_{n,i})\right]^\gamma \leq c_{2,2}\, n^6\, n^{-(1-\delta)\gamma} \to 0 \tag{28}$$

as $n \to \infty$. This implies that $\dim_{\rm H} U([0,1]^2) \leq \gamma$ a.s. By letting $\gamma_1 \uparrow 1/4$ and $\gamma_2 \uparrow 1/2$ along rational numbers, respectively, we derive (25).

Now we turn to the proof of the upper bound in (24) for the stationary pinned string $U$. We will show that there are three different ways to cover $\text{Gr} U([0,1]^2)$, each of which leads to an upper bound for $\dim_{\rm H} \text{Gr} U([0,1]^2)$.

- For each fixed integer $n \geq 2$, we have

$$\text{Gr} U([0,1]^2) \subseteq \bigcup_{i=1}^{n^6} R_{n,i} \times U(R_{n,i}). \tag{29}$$

It follows from (27) and (29) that $\text{Gr} U([0,1]^2)$ can be covered by $n^6$ cubes in $\mathbb{R}^{2+d}$ with side-lengths $c_{2,3}\, n^{-1+\delta}$ and the same argument as the above yields

$$\dim_{\rm H} \text{Gr} U([0,1]^2) \leq 6 \quad \text{a.s.} \tag{30}$$

- Observe that each $R_{n,i} \times U(R_{n,i})$ can be covered by $\ell_{n,1}$ cubes in $\mathbb{R}^{2+d}$ of sides $n^{-4}$, where by (26)

$$\ell_{n,1} \leq c_{2,4}\, n^2 \times \left(\frac{n^{-1+\delta}}{n^{-4}}\right)^d.$$

Hence $\text{Gr} U([0,1]^2)$ can be covered by $n^6 \times \ell_{n,1}$ cubes in $\mathbb{R}^{2+d}$ with sides $n^{-4}$. Denote

$$\eta_1 = 2 + (1-\gamma_1)d.$$

Recall from the above that we can choose the constants $\gamma_1$, $\gamma_1'$ and $\gamma_2'$ such that $1 - \delta > 4\gamma_1$. Therefore

$$n^6 \times \ell_{n,1} \times (n^{-4})^{\eta_1} \leq c_{2,5}\, n^{-(1-\delta-4\gamma_1)d} \to 0$$

as $n \to \infty$. This implies that $\dim_{\rm H} \text{Gr} U([0,1]^2) \leq \eta_1$ almost surely. Hence,

$$\dim_{\rm H} \text{Gr} U([0,1]^2) \leq 2 + \frac{3}{4}d, \quad \text{a.s.} \tag{31}$$

- We can also cover each $R_{n,i} \times U(R_{n,i})$ by $\ell_{n,2}$ cubes in $\mathbb{R}^{2+d}$ of sides $n^{-2}$, where by (26)

$$\ell_{n,2} \leq c_{2,6} \left(\frac{n^{-1+\delta}}{n^{-2}}\right)^d.$$

Hence $\mathrm{Gr}U([0,1]^2)$ can be covered by $n^6 \times \ell_{n,2}$ cubes in $\mathbb{R}^{2+d}$ with sides $n^{-2}$. Denote $\eta_2 = 3 + (1-\gamma_2)d$. Recall from the above that we can choose the constants $\gamma_2$, $\gamma_1'$ and $\gamma_2'$ such that $1 - \delta > 2\gamma_2$. Therefore

$$n^6 \times \ell_{n,2} \times (n^{-2})^{\eta_2} \leq c_{2,7} \, n^{-(1-\delta-2\gamma_2)d} \to 0$$

as $n \to \infty$. This implies that $\dim_{\mathrm{H}} \mathrm{Gr}U([0,1]^2) \leq \eta_2$ almost surely. Hence,

$$(32) \qquad \dim_{\mathrm{H}} \mathrm{Gr}U([0,1]^2) \leq 3 + \frac{1}{2}d, \qquad \text{a.s.}$$

Combining (30), (31) and (32) yields

$$(33) \qquad \dim_{\mathrm{H}} \mathrm{Gr}U([0,1]^2) \leq \min\left\{6, \, 2 + \frac{3}{4}d, \, 3 + \frac{1}{2}d\right\}, \qquad \text{a.s.}$$

and the upper bounds in (24) follow from (33).

To prove the lower bound in (23), by Frostman's theorem it is sufficient to show that for any $0 < \gamma < \min\{d, 6\}$,

$$(34) \qquad \mathcal{E}_\gamma = \int_{[0,1]^2} \int_{[0,1]^2} \mathbb{E}\left(\frac{1}{|U_s(y) - U_t(x)|^\gamma}\right) ds\, dy\, dt\, dx < \infty.$$

See, e.g., [7], Chapter 10. Since $0 < \gamma < d$, we have $0 < \mathbb{E}(|\Xi|^{-\gamma}) < \infty$, where $\Xi$ is a standard $d$-dimensional normal vector. Combining this fact with Lemma 1.1, we have

$$(35) \qquad \mathcal{E}_\gamma \leq c_{2,8} \int_0^1 ds \int_0^1 dt \int_0^1 dy \int_0^1 \frac{1}{\left(|s-t|^{1/2} + |x-y|\right)^{\gamma/2}} dx.$$

Recall the weighted arithmetic-mean and geometric-mean inequality: for all integer $n \geq 2$ and $x_i \geq 0$, $\beta_i > 0$ $(i = 1, \ldots, n)$ such that $\sum_{i=1}^n \beta_i = 1$, we have

$$(36) \qquad \prod_{i=1}^n x_i^{\beta_i} \leq \sum_{i=1}^n \beta_i x_i.$$

Applying (36) with $n = 2$, $\beta_1 = 2/3$ and $\beta_2 = 1/3$, we obtain

$$(37) \qquad |s-t|^{1/2} + |x-y| \geq \frac{2}{3}|s-t|^{1/2} + \frac{1}{3}|x-y| \geq |s-t|^{1/3}|x-y|^{1/3}.$$

Therefore, the denominator in (35) can be bounded from below by $|s-t|^{\gamma/6}|x-y|^{\gamma/6}$. Since $\gamma < 6$, by (35), we have $\mathcal{E}_\gamma < \infty$, which proves (34).

For proving the lower bound in (24), we need the following lemma from Ayache and Xiao [2].

**Lemma 2.2.** *Let $\alpha$, $\beta$ and $\eta$ be positive constants. For $a > 0$ and $b > 0$, let*

$$(38) \qquad J := J(a,b) = \int_0^1 \frac{dt}{(a + t^\alpha)^\beta (b+t)^\eta}.$$

*Then there exist finite constants $c_{2,9}$ and $c_{2,10}$, depending on $\alpha$, $\beta$, $\eta$ only, such that the following hold for all reals $a$, $b > 0$ satisfying $a^{1/\alpha} \leq c_{2,9} b$:*

(i) *if $\alpha\beta > 1$, then*

$$\text{(39)} \qquad J \leq c_{2,10} \, \frac{1}{a^{\beta-\alpha^{-1}} b^\eta};$$

(ii) *if $\alpha\beta = 1$, then*

$$\text{(40)} \qquad J \leq c_{2,10} \, \frac{1}{b^\eta} \log\left(1 + b a^{-1/\alpha}\right);$$

(iii) *if $0 < \alpha\beta < 1$ and $\alpha\beta + \eta \neq 1$, then*

$$\text{(41)} \qquad J \leq c_{2,10} \left( \frac{1}{b^{\alpha\beta+\eta-1}} + 1 \right).$$

Now we prove the lower bound in (24). Since $\dim_{\text{H}} \text{Gr} U([0,1]^2) \geq \dim_{\text{H}} U([0,1]^2)$ always holds, we only need to consider the cases $1 \leq d < 4$ and $4 \leq d < 6$, respectively.

Since the proof of the two cases are almost identical, we only prove the case when $1 \leq d < 4$ here. Let $0 < \gamma < 2 + \frac{3}{4}d$ be a fixed, but arbitrary, constant. Since $1 \leq d < 4$, we may and will assume $\gamma > 1+d$. In order to prove $\dim_{\text{H}} \text{Gr} U([0,1]^2) \geq \gamma$ a.s., again by Frostman's theorem, it is sufficient to show

$$\text{(42)} \qquad \mathcal{G}_\gamma = \int_{[0,1]^2} \int_{[0,1]^2} \mathbb{E}\left[ \frac{1}{\left(|s-t|^2 + |x-y|^2 + |U_s(y) - U_t(x)|^2\right)^{\gamma/2}} \right] ds\,dy\,dt\,dx$$
$$< \infty.$$

Since $\gamma > d$, we note that for a standard normal vector $\Xi$ in $\mathbb{R}^d$ and any number $a \in \mathbb{R}$,

$$\mathbb{E}\left[ \frac{1}{\left(a^2 + |\Xi|^2\right)^{\gamma/2}} \right] \leq c_{2,11} \, a^{-(\gamma-d)},$$

see e.g. [7], p. 279. Consequently, by Lemma 1.1, we derive that

$$\text{(43)} \quad \mathcal{G}_\gamma \leq c_{2,12} \int_{[0,1]^2} \int_{[0,1]^2} \frac{1}{\left(|s-t|^{1/2} + |x-y|\right)^{d/2} \left(|s-t| + |x-y|\right)^{\gamma-d}} ds\,dy\,dt\,dx.$$

By Lemma 2.2 and a change of variable and noting that $d < 4$, we can apply (41) to derive

$$\text{(44)} \qquad \begin{aligned} \mathcal{G}_\gamma &\leq c_{2,13} \int_0^1 dx \int_0^1 \frac{1}{(t^{1/2} + x)^{d/2} (t+x)^{\gamma-d}} dt \\ &\leq c_{2,14} \int_0^1 \left( \frac{1}{x^{d/4+\gamma-d-1}} + 1 \right) dx < \infty, \end{aligned}$$

where the last inequality follows from $\gamma - \frac{3}{4}d - 1 < 1$. This completes the proof of Theorem 2.1. $\square$

By using the relationships among the Hausdorff dimension, packing dimension and the box dimension (see Falconer [4]), Theorem 2.1 and the proof of the upper bounds, we derive the following analogous result on the packing dimensions of $u([0,1]^2)$ and $\text{Gr} u([0,1]^2)$.

**Theorem 2.3.** *Let $\{u_t(x) : t \geq 0, x \in \mathbb{R}\}$ be a random string process taking values in $\mathbb{R}^d$. Then with probability 1,*

$$\dim_{\mathrm{P}} u([0,1]^2) = \min\{d; 6\} \tag{45}$$

*and*

$$\dim_{\mathrm{P}} \mathrm{Gr} u([0,1]^2) = \begin{cases} 2 + \frac{3}{4}d & \text{if } 1 \leq d < 4, \\ 3 + \frac{1}{2}d & \text{if } 4 \leq d < 6, \\ 6 & \text{if } 6 \leq d. \end{cases} \tag{46}$$

Theorems 2.1 and 2.3 show that the random fractals $u([0,1]^2)$ and $\mathrm{Gr} u([0,1]^2)$ are rather regular because they have the same Hausdorff and packing dimensions.

Now we will turn our attention to find the Hausdorff dimension of the range $u(E)$ for an arbitrary Borel set $E \subseteq [0,\infty) \times \mathbb{R}$.

For this purpose, we mention the related results of Wu and Xiao [13] for an $(N,d)$-fractional Brownian sheet $B^H = \{B^H(\mathbf{t}) : \mathbf{t} \in \mathbb{R}_+^N\}$ with Hurst index $H = (H_1, \ldots, H_N) \in (0,1)^N$. What the random string process $\{u_t(x) : t \geq 0, x \in \mathbb{R}\}$ and a $(2,d)$-fractional Broanian sheet $B^H$ with $H = (\frac{1}{4}, \frac{1}{2})$ have in common is that they are both anisotropic.

As Wu and Xiao [13] pointed out, the Hausdorff dimension of the image $B^H(F)$ cannot be determined by $\dim_{\mathrm{H}} F$ and $H$ alone for an arbitrary fractal set $F$, and more information about the geometry of $F$ is needed. To capture the anisotropic nature of $B^H$, they have introduced a new notion of dimension, namely, the *Hausdorff dimension contour*, for finite Borel measures and Borel sets and showed that $\dim_{\mathrm{H}} B^H(F)$ is determined by the Hausdorff dimension contour of $F$. It turns out that we can use the same technique to study the images of the random string.

We start with the following Proposition 2.4 which determines $\dim_{\mathrm{H}} u(E)$ when $E$ belongs to a special class of Borel sets in $[0,\infty) \times \mathbb{R}$. Its proof is the same as that of Proposition 3.1 in [13].

**Proposition 2.4.** *Let $\{u_t(x) : t \geq 0, x \in \mathbb{R}\}$ be a random string in $\mathbb{R}^d$. Assume that $E_1$ and $E_2$ are Borel sets in $[0,\infty)$ and $\mathbb{R}$, respectively, which satisfy $\dim_{\mathrm{H}} E_1 = \dim_{\mathrm{P}} E_1$ or $\dim_{\mathrm{H}} E_2 = \dim_{\mathrm{P}} E_2$. Let $E = E_1 \times E_2 \subset [0,\infty) \times \mathbb{R}$, then we have*

$$\dim_{\mathrm{H}} u(E) = \min\{d; 4\dim_{\mathrm{H}} E_1 + 2\dim_{\mathrm{H}} E_2\}, \quad a.s. \tag{47}$$

In order to determine $\dim_{\mathrm{H}} u(E)$ for an arbitrary Borel set $E \subset [0,\infty) \times \mathbb{R}$, we recall from [13] the following definition. Denote by $\mathcal{M}_c^+(E)$ the family of finite Borel measures with compact support in $E$.

**Definition 2.5.** *Given $\mu \in \mathcal{M}_c^+(E)$, we define the set $\Lambda_\mu \subseteq \mathbb{R}_+^2$ by*

$$\Lambda_\mu = \left\{ \lambda = (\lambda_1, \lambda_2) \in \mathbb{R}_+^2 : \limsup_{r \to 0^+} \frac{\mu(R((t,x),r))}{r^{4\lambda_1 + 2\lambda_2}} = 0, \right. \tag{48}$$
$$\left. \text{for } \mu\text{-a.e. } (t,x) \in [0,\infty) \times \mathbb{R} \right\},$$

*where $R((t,x), r) = [t - r^4, t + r^4] \times [x - r^2, x + r^2]$.*

The properties of set $\Lambda_\mu$ can be found in Lemma 3.6 of Wu and Xiao [13]. The boundary of $\Lambda_\mu$, denoted by $\partial \Lambda_\mu$, is called the Hausdorff dimension contour of $\mu$.

Define
$$\Lambda(E) = \bigcup_{\mu \in \mathcal{M}_c^+(E)} \Lambda_\mu.$$

and define the Hausdorff dimension contour of $E$ by $\bigcup_{\mu \in \mathcal{M}_c^+(E)} \partial \Lambda_\mu$. It can be verified that, for every $\mathbf{b} \in (0, \infty)^2$, the supremum $\sup_{\lambda \in \Lambda(E)} \langle \lambda, \mathbf{b} \rangle$ is achieved on the Hausdorff dimension contour of $E$ (Lemma 3.6, [13]).

**Theorem 2.6.** *Let $u = \{u_t(x) : t \geq 0, x \in \mathbb{R}\}$ be a random string process with values in $\mathbb{R}^d$. Then, for any Borel set $E \subset [0, \infty) \times \mathbb{R}$,*

$$\dim_{\mathrm{H}} u(E) = \min\{d; s(E)\}, \qquad a.s. \tag{49}$$

*where $s(E) = \sup_{\lambda \in \Lambda(E)}(4\lambda_1 + 2\lambda_2)$.*

*Proof.* By Corollary 2 of Mueller and Tribe [9], one only needs to prove (49) for the stationary pinned string $U = \{U_t(x) : t \geq 0, x \in \mathbb{R}\}$. The latter follows from the proof of Theorem 3.10 in [13]. □

## 3. Existence of the local times and dimension results for level sets

In this section, we will first give a sufficient condition for the existence of the local times of a random string process on any rectangle $I \in \mathcal{A}$, where $\mathcal{A}$ is the collection of all the rectangles in $[0, \infty) \times \mathbb{R}$ with sides parallel to the axes. Then, we will determine the Hausdorff and packing dimensions for the level set $L_{\mathbf{u}} = \{(t, x) \in [0, \infty) \times \mathbb{R} : u_t(x) = \mathbf{u}\}$, where $\mathbf{u} \in \mathbb{R}^d$ is fixed.

We start by briefly recalling some aspects of the theory of local times. For an excellent survey on local times of random and deterministic vector fields, we refer to Geman and Horowitz [6].

Let $X(\mathbf{t})$ be a Borel vector field on $\mathbb{R}^N$ with values in $\mathbb{R}^d$. For any Borel set $T \subseteq \mathbb{R}^N$, the occupation measure of $X$ on $T$ is defined as the following measure on $\mathbb{R}^d$:
$$\mu_T(\bullet) = \lambda_N\{\mathbf{t} \in T : X(\mathbf{t}) \in \bullet\}.$$

If $\mu_T$ is absolutely continuous with respect to $\lambda_d$, the Lebesgue measure on $\mathbb{R}^d$, we say that $X(\mathbf{t})$ has *local times* on $T$, and define its local time $l(\bullet, T)$ as the Radon–Nikodým derivative of $\mu_T$ with respect to $\lambda_d$, i.e.,

$$l(\mathbf{u}, T) = \frac{d\mu_T}{d\lambda_d}(\mathbf{u}), \qquad \forall \mathbf{u} \in \mathbb{R}^d.$$

In the above, $\mathbf{u}$ is the so-called *space variable*, and $T$ is the *time* variable. Sometimes, we write $l(\mathbf{u}, \mathbf{t})$ in place of $l(\mathbf{u}, [0, \mathbf{t}])$. It is clear that if $X$ has local times on $T$, then for every Borel set $S \subseteq T$, $l(\mathbf{u}, S)$ also exists.

By standard martingale and monotone class arguments, one can deduce that the local times have a measurable modification that satisfies the following *occupation density formula*: for every Borel set $T \subseteq \mathbb{R}^N$, and for every measurable function $f : \mathbb{R}^d \to \mathbb{R}$,

$$\int_T f(X(\mathbf{t})) \, d\mathbf{t} = \int_{\mathbb{R}^d} f(\mathbf{u}) l(\mathbf{u}, T) \, d\mathbf{u}. \tag{50}$$

The following theorem is concerned with the existence of local times of the random string.

**Theorem 3.1.** *Let $\{u_t(x) : t \geq 0, x \in \mathbb{R}\}$ be a random string process in $\mathbb{R}^d$. If $d < 6$, then for every $I \in \mathcal{A}$, the string has local times $\{l(\mathbf{u}, I), \mathbf{u} \in \mathbb{R}^d\}$ on $I$, and $l(\mathbf{u}, I)$ admits the following $L^2$ representation:*

$$(51) \qquad l(\mathbf{u}, I) = (2\pi)^{-d} \int_{\mathbb{R}^d} e^{-i\langle \mathbf{v}, \mathbf{u} \rangle} \int_I e^{i\langle \mathbf{v}, u_t(x) \rangle} dt dx d\mathbf{v}, \quad \forall\, \mathbf{u} \in \mathbb{R}^d.$$

*Proof.* Because of Corollary 2 of Mueller and Tribe [9], we only need to prove the existence for the stationary pinned string $U = \{U_t(x) : t \geq 0, x \in \mathbb{R}\}$.

Let $I \in \mathcal{A}$ be fixed. Without loss of generality, we may assume $I = [\varepsilon, 1]^2$. By (21.3) in [6] and using the characteristic functions of Gaussian random variables, it suffices to prove

$$(52) \quad \mathcal{J}(I) := \int_I dt dx \int_I ds dy \int_{\mathbb{R}^d} d\mathbf{u} \int_{\mathbb{R}^d} \left| \mathbb{E} \exp\left(i\langle \mathbf{u}, U_t(x) \rangle + i\langle \mathbf{v}, U_s(y) \rangle\right) \right| d\mathbf{v} < \infty.$$

Since the components of $U$ are i.i.d., it is easy to see that

$$(53) \qquad \mathcal{J}(I) = (2\pi)^d \int_I dt dx \int_I \left[\det\text{Cov}\left(U_t^1(x), U_s^1(y)\right)\right]^{-d/2} ds dy.$$

By Lemma 1.3 and noting that $I = [\varepsilon, 1]^2$, we can see that

$$(54) \qquad \begin{aligned} \det\text{Cov}\left(U_t^j(x), U_s^j(y)\right) &= \text{Var}\left(U_s^j(y)\right)\text{Var}\left(U_t^j(x)|U_s^j(y)\right) \\ &\geq c_{3,1}\left(|x-y| + |t-s|^{1/2}\right). \end{aligned}$$

The above inequality, (37) and the fact that $d < 6$ lead to

$$(55) \qquad \mathcal{J}(I) \leq c_{3,2} \int_\varepsilon^1 \int_\varepsilon^1 |s-t|^{-d/6} dt ds \int_\varepsilon^1 \int_\varepsilon^1 |x-y|^{-d/6} dx dy < \infty,$$

which proves (52), and therefore Theorem 3.1. $\square$

**Remark 3.2.** It would be interesting to study the regularity properties of the local times $l(\mathbf{u}, \mathbf{t})$, $(\mathbf{u} \in \mathbb{R}^d, \mathbf{t} \in [0, \infty) \times \mathbb{R})$ such as joint continuity and moduli of continuity. One way to tackle these problems is to establish sectorial local nondeterminism (see [13]) for the stationary pinned string $U = \{U_t(x) : t \geq 0, x \in \mathbb{R}\}$. This will have to be pursued elsewhere. Some results of this nature for certain isotropic Gaussian random fields can be found in [15].

Mueller and Tribe [9] proved that for every $\mathbf{u} \in \mathbb{R}^d$,

$$\mathbb{P}\{u_t(x) = \mathbf{u} \text{ for some } (t, x) \in [0, \infty) \times \mathbb{R}\} > 0$$

if and only if $d < 6$. Now we study the Hausdorff and packing dimensions of the level set $L_\mathbf{u} = \{(t, x) \in [0, \infty) \times \mathbb{R} : u_t(x) = \mathbf{u}\}$.

**Theorem 3.3.** *Let $\{u_t(x) : t \geq 0, x \in \mathbb{R}\}$ be a random string process in $\mathbb{R}^d$ with $d < 6$. Then for every $\mathbf{u} \in \mathbb{R}^d$, with positive probability,*

$$(56) \qquad \dim_{\text{H}}\left(L_\mathbf{u} \cap [0,1]^2\right) = \dim_{\text{P}}\left(L_\mathbf{u} \cap [0,1]^2\right) = \begin{cases} 2 - \frac{1}{4}d & \text{if } 1 \leq d < 4, \\ 3 - \frac{1}{2}d & \text{if } 4 \leq d < 6. \end{cases}$$

*Proof.* As usual, it is sufficient to prove (56) for the stationary pinned string $U = \{U_t(x) : t \geq 0, \ x \in \mathbb{R}\}$. We first prove the almost sure upper bound

(57) $$\dim_{\mathrm{P}}\bigl(L_{\mathbf{u}} \cap [0,1]^2\bigr) \leq \begin{cases} 2 - \frac{1}{4}d & \text{if } 1 \leq d < 4, \\ 3 - \frac{1}{2}d & \text{if } 4 \leq d < 6. \end{cases}$$

By the $\sigma$-stability of $\dim_{\mathrm{P}}$, it is sufficient to show (57) holds for $L_{\mathbf{u}} \cap [\varepsilon, 1]^2$ for every $\varepsilon \in (0,1)$. For this purpose, we construct coverings of $L_{\mathbf{u}} \cap [0,1]^2$ by cubes of the same side length.

For any integer $n \geq 2$, we divide the square $[\varepsilon, 1]^2$ into $n^6$ sub-rectangles $R_{n,\ell}$ of side lengths $n^{-4}$ and $n^{-2}$, respectively. Let $0 < \delta < 1$ be fixed and let $\tau_{n,\ell}$ be the lower-left vertex of $R_{n,\ell}$. Then

(58)
$$\begin{aligned}
\mathbb{P}\bigl\{\mathbf{u} \in U(R_{n,\ell})\bigr\} &\leq \mathbb{P}\Bigl\{\max_{(s,y),(t,x)\in R_{n,\ell}} |U_s(y) - U_t(x)| \leq n^{-(1-\delta)};\ \mathbf{u} \in U(R_{n,\ell})\Bigr\} \\
&\quad + \mathbb{P}\Bigl\{\max_{(s,y),(t,x)\in R_{n,\ell}} |U_s(y) - U_t(x)| > n^{-(1-\delta)}\Bigr\} \\
&\leq \mathbb{P}\Bigl\{|U(\tau_{n,\ell}) - \mathbf{u}| \leq n^{-(1-\delta)}\Bigr\} + e^{-c n^{2\delta}} \\
&\leq c_{3,3}\, n^{-(1-\delta)d}.
\end{aligned}$$

In the above we have applied Lemma 1.1 and the Gaussian isoperimetric inequality (cf. Lemma 2.1 in [11]) to derive the second inequality.

Since we can deal with the cases $1 \leq d < 4$ and $4 \leq d < 6$ almost identically, we will only consider the case $1 \leq d < 4$ here and leave the case $4 \leq d < 6$ to the interested readers.

Define a covering $\{R'_{n,\ell}\}$ of $L_{\mathbf{u}} \cap [\varepsilon, 1]^2$ by $R'_{n,\ell} = R_{n,\ell}$ if $\mathbf{u} \in U(R_{n,\ell})$ and $R'_{n,\ell} = \emptyset$ otherwise. Note that each $R'_{n,\ell}$ can be covered by $n^2$ squares of side length $n^{-4}$. Thus, for every $n \geq 2$, we have obtained a covering of the level set $L_{\mathbf{u}} \cap [\varepsilon, 1]^2$ by squares of side length $n^{-4}$. Consider the sequence of integers $n = 2^k$ ($k \geq 1$), and let $N_k$ denote the minimum number of squares of side-length $2^{-4k}$ that are needed to cover $L_{\mathbf{u}} \cap [\varepsilon, 1]^2$. It follows from (58) that

(59) $$\mathbb{E}(N_k) \leq c_{3,3}\, 2^{6k} \cdot 2^{2k} \cdot 2^{-k(1-\delta)d} = c_{3,3}\, 2^{k(8-(1-\delta)d)}.$$

By (59), Markov's inequality and the Bore-Cantelli lemma we derive that for any $\delta' \in (0, \delta)$, almost surely for all $k$ large enough,

(60) $$N_k \leq c_{3,3}\, 2^{k(8-(1-\delta')d)}.$$

By the definition of box dimension and its relation to $\dim_{\mathrm{P}}$ (cf. [4]), (60) implies that $\dim_{\mathrm{P}}\bigl(L_{\mathbf{u}} \cap [\varepsilon, 1]^2\bigr) \leq 2 - (1-\delta')d/4$ a.s. Since $\varepsilon > 0$ is arbitrary, we obtain the desired upper bound for $\dim_{\mathrm{P}}\bigl(L_{\mathbf{u}} \cap [\varepsilon, 1]^2\bigr)$ in the case $1 \leq d < 4$.

Since $\dim_{\mathrm{H}} E \leq \dim_{\mathrm{P}} E$ for all Borel sets $E \subset \mathbb{R}^2$, it remains to prove the following lower bound: for any $\varepsilon \in (0,1)$, with positive probability

(61) $$\dim_{\mathrm{P}}\bigl(L_{\mathbf{u}} \cap [\varepsilon, 1]^2\bigr) \geq \begin{cases} 2 - \frac{1}{4}d & \text{if } 1 \leq d < 4, \\ 3 - \frac{1}{2}d & \text{if } 4 \leq d < 6. \end{cases}$$

We only prove (61) for $1 \leq d < 4$. The other case is similar and is omitted. Let $\delta > 0$ such that

(62) $$\gamma := 2 - \frac{1}{4}(1+\delta)d > 1.$$

Note that if we can prove that there is a constant $c_{3,4} > 0$ such that

(63) $$\mathbb{P}\{\dim_{\mathrm{H}}(L_{\mathbf{u}} \cap [\varepsilon, 1]^2) \geq \gamma\} \geq c_{3,4},$$

then the lower bound in (61) will follow by letting $\delta \downarrow 0$.

Our proof of (63) is based on the capacity argument due to Kahane (see, e.g., [7]). Similar methods have been used by Adler [1], Testard [12], Xiao [14], Ayache and Xiao [2] to various types of stochastic processes.

Let $\mathcal{M}_\gamma^+$ be the space of all non-negative measures on $[0,1]^2$ with finite $\gamma$-energy. It is known (cf. [1]) that $\mathcal{M}_\gamma^+$ is a complete metric space under the metric

(64) $$\|\mu\|_\gamma = \int_{\mathbb{R}^2} \int_{\mathbb{R}^2} \frac{\mu(dt, dx)\mu(ds, dy)}{(|t-s|^2 + |x-y|^2)^{\gamma/2}}.$$

We define a sequence of random positive measures $\mu_n$ on the Borel sets of $[\varepsilon, 1]^2$ by

(65) $$\begin{aligned} \mu_n(C) &= \int_C (2\pi n)^{d/2} \exp\left(-\frac{n|U_t(x) - \mathbf{u}|^2}{2}\right) dt dx \\ &= \int_C \int_{\mathbb{R}^d} \exp\left(-\frac{|\xi|^2}{2n} + i\langle \xi, U_t(x) - \mathbf{u}\rangle\right) d\xi\, dtdx, \ \forall\, C \in \mathcal{B}([\varepsilon,1]^2). \end{aligned}$$

It follows from Kahane [7] or Testard [12] that if there are positive constants $c_{3,5}$ and $c_{3,6}$, which may depend on $\mathbf{u}$, such that

(66) $$\mathbb{E}(\|\mu_n\|) \geq c_{3,5}, \qquad \mathbb{E}(\|\mu_n\|^2) \leq c_{3,6},$$
(67) $$\mathbb{E}(\|\mu_n\|_\gamma) < +\infty,$$

where $\|\mu_n\| = \mu_n([\varepsilon, 1]^2)$, then there is a subsequence of $\{\mu_n\}$, say $\{\mu_{n_k}\}$, such that $\mu_{n_k} \to \mu$ in $\mathcal{M}_\gamma^+$ and $\mu$ is strictly positive with probability $\geq c_{3,5}^2/(2c_{3,6})$. It follows from (65) and the continuity of $U$ that $\mu$ has its support in $L_\mathbf{u} \cap [\varepsilon, 1]^2$ almost surely. Hence Frostman's theorem yields (63).

It remains to verify (66) and (67). By Fubini's theorem we have

(68) $$\begin{aligned} \mathbb{E}(\|\mu_n\|) &= \int_{[\varepsilon,1]^2} \int_{\mathbb{R}^d} e^{-i\langle\xi,\mathbf{u}\rangle} \exp\left(-\frac{|\xi|^2}{2n}\right) \mathbb{E}\exp\left(i\langle\xi, U_t(x)\rangle\right) d\xi\, dtdx \\ &= \int_{[\varepsilon,1]^2} \int_{\mathbb{R}^d} e^{-i\langle\xi,\mathbf{u}\rangle} \exp\left(-\frac{1}{2}(n^{-1} + \sigma^2(t,x))|\xi|^2\right) d\xi\, dtdx \\ &= \int_{[\varepsilon,1]^2} \left(\frac{2\pi}{n^{-1} + \sigma^2(t,x)}\right)^{d/2} \exp\left(-\frac{|\mathbf{u}|^2}{2(n^{-1}+\sigma^2(t,x))}\right) dtdx \\ &\geq \int_{[\varepsilon,1]^2} \left(\frac{2\pi}{1+\sigma^2(t,x)}\right)^{d/2} \exp\left(-\frac{|\mathbf{u}|^2}{2\sigma^2(t,x)}\right) dt := c_{3,5}, \end{aligned}$$

where $\sigma^2(t,x) = \mathbb{E}\left[(U_t^1(x))^2\right]$.

Denote by $I_{2d}$ the identity matrix of order $2d$ and by $\mathrm{Cov}(U_s(y), U_t(x))$ the covariance matrix of the Gaussian vector $(U_s(y), U_t(x))$. Let $\Gamma = n^{-1}I_{2d} + \mathrm{Cov}(U_s(y), U_t(x))$ and let $(\xi,\eta)'$ be the transpose of the row vector $(\xi,\eta)$. As in the proof of

(52), we apply (14) in Lemma 1.3 and the inequality (36) to derive

$$
\begin{aligned}
&\mathbb{E}\big(\|\mu_n\|^2\big) \\
&= \int_{[\varepsilon,1]^2}\int_{[\varepsilon,1]^2}\int_{\mathbb{R}^d}\int_{\mathbb{R}^d} e^{-i\langle\xi+\eta,\mathbf{u}\rangle}\exp\Big(-\frac{1}{2}(\xi,\eta)\,\Gamma\,(\xi,\eta)'\Big)\,d\xi d\eta\,dsdydtdx \\
&= \int_{[\varepsilon,1]^2}\int_{[\varepsilon,1]^2}\frac{(2\pi)^d}{\sqrt{\det\Gamma}}\exp\Big(-\frac{1}{2}(\mathbf{u},\mathbf{u})\,\Gamma^{-1}\,(\mathbf{u},\mathbf{u})'\Big)\,dsdy\,dtdx \\
&\leq \int_{[\varepsilon,1]^N}\int_{[\varepsilon,1]^N}\frac{(2\pi)^d}{\big[\det\mathrm{Cov}(U_s^1(y),U_t^1(x))\big]^{d/2}}\,dsdy\,dtdx \\
&\leq c_{3,7}\int_\varepsilon^1\int_\varepsilon^1 |s-t|^{-d/6}dtds\int_\varepsilon^1\int_\varepsilon^1 |x-y|^{-d/6}dxdy := c_{3,6} < \infty.
\end{aligned} \tag{69}
$$

Similar to (69) and by the same method as in proving (43), we have

$$
\begin{aligned}
\mathbb{E}(\|\mu_n\|_\gamma) &= \int_{[\varepsilon,1]^2}\int_{[\varepsilon,1]^2}\frac{dsdy\,dtdx}{(|s-t|^2+|x-y|^2)^{\gamma/2}} \\
&\quad \times \int_{\mathbb{R}^d}\int_{\mathbb{R}^d} e^{-i\langle\xi+\eta,\mathbf{u}\rangle}\exp\Big(-\frac{1}{2}(\xi,\eta)\,\Gamma\,(\xi,\eta)'\Big)\,d\xi d\eta \\
&\leq c_{3,8}\int_{[\varepsilon,1]^2}\int_{[\varepsilon,1]^2}\frac{dsdy\,dtdx}{\big(|s-t|^{1/2}+|x-y|\big)^{d/2}\big(|s-t|+|x-y|\big)^\gamma} \\
&< \infty,
\end{aligned} \tag{70}
$$

where the last inequality follows from Lemma 2.2 and the facts that $d < 4$ and $d/4+\gamma-1 < 1$. This proves (67) and thus the proof of Theorem 3.3 is finished. □

## 4. Hausdorff and packing dimensions of the sets of double times

Mueller and Tribe [9] found necessary and sufficient conditions for an $\mathbb{R}^d$-valued string process to have double points. In this section, we determine the Hausdorff and packing dimensions of the sets of double times of the random string.

As in [9], we consider the following two kinds of double times for the string process $\{u_t(x) : t \geq 0,\, x \in \mathbb{R}\}$.

- Type I double times:

$$
L_{\mathrm{I},2} = \Big\{\big((t_1,x_1),(t_2,x_2)\big) \in \big((0,\infty)\times\mathbb{R}\big)_{\neq}^2 :\ u_{t_1}(x_1) = u_{t_2}(x_2)\Big\}, \tag{71}
$$

where

$$
\big((0,\infty)\times\mathbb{R}\big)_{\neq}^2 = \big\{\big((t_1,x_1),(t_2,x_2)\big) \in \big((0,\infty)\times\mathbb{R}\big)^2 :\ (t_1,x_1) \neq (t_2,x_2)\big\}.
$$

In order to determine the Hausdorff and packing dimensions of $L_{\mathrm{I},2}$, we introduce a $(4,d)$-random field $\Delta u = \{\Delta u(t_1,x_1;t_2,x_2)\}$ defined by

$$
\Delta u(t_1,x_1;t_2,x_2) = u_{t_2}(x_2) - u_{t_1}(x_1),\ \forall (t_1,x_1,t_2,x_2) \in \big((0,\infty)\times\mathbb{R}\big)^2. \tag{72}
$$

Then $L_{\mathrm{I},2}$ can be viewed as the zero set of $\Delta u(t_1,x_1;t_2,x_2)$, denoted by $(\Delta u)^{-1}(0)$; and its Hausdorff and packing dimensions can be studied by using the method in Section 3.

- Type II double times:

$$L_{\text{II},2} = \left\{(t, x_1, x_2) \in (0, \infty) \times \mathbb{R}^2_{\neq} : u_t(x_1) = u_t(x_2)\right\},$$

where $\mathbb{R}^2_{\neq} = \{(x_1, x_2) \in \mathbb{R}^2 : x_1 \neq x_2\}$.

In order to determine the Hausdorff and packing dimensions of $L_{\text{II},2}$, we will consider the $(3, d)$-random field $\widetilde{\Delta} u = \{\widetilde{\Delta} u(t; x_1, x_2)\}$ defined by

(73) $\quad \widetilde{\Delta} u(t; x_1, x_2) = u_t(x_2) - u_t(x_1), \quad \forall \, (t, x_1, x_2) \in (0, \infty) \times \mathbb{R}^2.$

Then we can see that $L_{\text{II},2}$ is nothing but the zero set of $\widetilde{\Delta} u$:

$$L_{\text{II},2} = \left\{(t, x_1, x_2) \in (0, \infty) \times \mathbb{R}^2_{\neq} : \widetilde{\Delta} u(t; x_1, x_2) = 0\right\}.$$

For any constants $0 < a_1 < a_2$ and $b_1 < b_2$, consider the squares $J_\ell = [a_\ell, a_\ell + h] \times [b_\ell, b_\ell + h]$ ($\ell = 1, 2$). Let $J = \prod_{\ell=1}^{2} J_\ell \subset ((0, \infty) \times \mathbb{R})^2$ denote the corresponding hypercube. We choose $h > 0$ small enough, say,

$$h < \min\left\{\frac{a_2 - a_1}{3}, \frac{b_2 - b_1}{3}\right\} \equiv L.$$

Thus $|t_2 - t_1| > L$ for all $t_2 \in [a_2, a_2 + h]$ and $t_1 \in [a_1, a_1 + h]$. We will use this assumption together with Lemma 1.4 to prove Theorem 4.1 below. We denote the collection of the hypercubes having the above properties by $\mathcal{J}$.

The following theorem gives the Hausdorff and packing dimensions of the Type I double times of a random string.

**Theorem 4.1.** *Let $u = \{u_t(x) : t \geq 0, x \in \mathbb{R}\}$ be a random string process in $\mathbb{R}^d$. If $d \geq 12$, then $L_{\text{I},2} = \emptyset$ a.s. If $d < 12$, then, for every $J \in \mathcal{J}$, with positive probability,*

(74) $\quad \dim_{\text{H}}\left(L_{\text{I},2} \cap J\right) = \dim_{\text{P}}\left(L_{\text{I},2} \cap J\right) = \begin{cases} 4 - \frac{1}{4}d & \text{if } 1 \leq d < 8, \\ 6 - \frac{1}{2}d & \text{if } 8 \leq d < 12. \end{cases}$

*Proof.* The first statement is due to Mueller and Tribe [9]. Hence, we only need to prove the dimension result (74).

Thanks to Corollary 2 of Mueller and Tribe [9], it is sufficient to prove (74) for the stationary pinned string $U$. This will be done by working with the zero set of the $(4, d)$-Gaussian field $\Delta U = \{\Delta U(t_1, x_1; t_2, x_2)\}$ define by (72). That is, we will prove (74) with $L_{\text{I},2}$ replaced by the zero set $(\Delta U)^{-1}(0)$. The proof is a modification of that of Theorem 3.3. Hence, we only give a sketch of it.

For an integer $n \geq 2$, we divide the hypercube $J$ into $n^{12}$ sub-domains $T_{n,p} = R^1_{n,p} \times R^2_{n,p}$, where $R^1_{n,p}, R^2_{n,p} \subset (0, \infty) \times \mathbb{R}$ are rectangles of side lengths $n^{-4}h$ and $n^{-2}h$, respectively. Let $0 < \delta < 1$ be fixed and let $\tau^k_{n,p}$ be the lower-left vertex of $R^k_{n,p}$ ($k = 1, 2$). Denote

$$\Delta V^{t_1, x_1; t_2, x_2}_{s_1, y_1; s_2, y_2} = \Delta U(t_1, x_1; t_2, x_2) - \Delta U(s_1, y_1; s_2, y_2),$$

then the probability $\mathbb{P}\{0 \in \Delta U(T_{n,p})\}$ is at most

$$
\begin{aligned}
&\mathbb{P}\bigg\{\max_{(t_1,x_1;t_2,x_2),(s_1,y_1;s_2,y_2)\in T_{n,p}}\big|\Delta V^{t_1,x_1;t_2,x_2}_{s_1,y_1;s_2,y_2}\big| \leq n^{-(1-\delta)};\, 0 \in \Delta U(T_{n,p})\bigg\}\\
&\quad + \mathbb{P}\bigg\{\max_{(t_1,x_1;t_2,x_2),(s_1,y_1;s_2,y_2)\in T_{n,p}}\big|\Delta V^{t_1,x_1;t_2,x_2}_{s_1,y_1;s_2,y_2}\big| > n^{-(1-\delta)}\bigg\}\\
&\leq \mathbb{P}\bigg\{\big|\Delta U(\tau^1_{n,p};\tau^2_{n,p})\big| \leq n^{-(1-\delta)}\bigg\}\\
&\quad + \mathbb{P}\bigg\{\max_{(t_1,x_1;t_2,x_2),(s_1,y_1;s_2,y_2)\in T_{n,p}}\big|\Delta V^{t_1,x_1;t_2,x_2}_{s_1,y_1;s_2,y_2}\big| > n^{-(1-\delta)}\bigg\}.
\end{aligned}
\tag{75}
$$

By the definition of $J$, we see that $\Delta U(\tau^1_{n,p}, \tau^2_{n,p})$ is a Gaussian random variable with mean 0 and variance at least $c\,L^{1/2}$. Hence the first term in (75) is at most $c_{4,1}\, n^{-(1-\delta)d}$.

On the other hand, since

$$\big|\Delta V^{t_1,x_1;t_2,x_2}_{s_1,y_1;s_2,y_2}\big| \leq c \sum_{k=1}^{2} \big|U_{s_k}(y_k) - U_{t_k}(x_k)\big|,$$

we have

$$
\begin{aligned}
&\mathbb{P}\bigg\{\max_{(t_1,x_1;t_2,x_2),(s_1,y_1;s_2,y_2)\in T_{n,p}}\big|\Delta V^{t_1,x_1;t_2,x_2}_{s_1,y_1;s_2,y_2}\big| > n^{-(1-\delta)}\bigg\}\\
&\leq \sum_{k=1}^{2}\mathbb{P}\bigg\{\max_{(s_k,y_k),(t_k,x_k)\in R^k_{n,p}}\big|U_{s_k}(y_k) - U_{t^k}(x_k)\big| > \frac{n^{-(1-\delta)}}{2c}\bigg\}\\
&\leq e^{-c_{4,2}\, n^{2\delta}},
\end{aligned}
\tag{76}
$$

where the last inequality follows from Lemma 1.1 and the Gaussian isoperimetric inequality (cf. Lemma 2.1 in [11]).

Combine (75) and (76), we have

$$\mathbb{P}\big\{0 \in \Delta U(T_{n,p})\big\} \leq c_{4,1}\, n^{-(1-\delta)d} + e^{-c_{4,2}\, n^{2\delta}}. \tag{77}$$

Hence the same covering argument as in the proof of Theorem 3.3 yields the desired upper bound for $\dim_{\mathrm{P}}\big((\Delta U)^{-1}(0)\cap J\big)$. This proves the upper bounds in (74).

Now we prove the lower bound for the Hausdorff dimension of $(\Delta U)^{-1}(0)\cap J$. We will only consider the case $1 \leq d < 8$ here and leave the case $8 \leq d < 12$ to the interested readers.

Let $\delta > 0$ such that

$$\gamma := 4 - \frac{1}{4}(1+\delta)d > 2. \tag{78}$$

As in the proof of Theorem 3.3, it is sufficient to prove that there is a constant $c_{4,3} > 0$ such that

$$\mathbb{P}\big\{\dim_{\mathrm{H}}\big(L_{\mathrm{I},2}\cap J\big) \geq \gamma\big\} \geq c_{4,3}. \tag{79}$$

Let $\mathcal{N}^+_\gamma$ be the space of all non-negative measures on $[0,1]^4$ with finite $\gamma$-energy. Then $\mathcal{N}^+_\gamma$ is a complete metric space under the metric

$$\|\nu\|_\gamma = \int_{\mathbb{R}^4}\int_{\mathbb{R}^4} \frac{\nu(dt_1 dx_1 dt_2 dx_2)\nu(ds_1 dy_1 ds_2 dy_2)}{(|t_1-s_1|^2 + |x_1-y_1|^2 + |t_2-s_2|^2 + |x_2-y_2|^2)^{\gamma/2}}; \tag{80}$$

see [1]. We define a sequence of random positive measures $\nu_n$ on the Borel set $J$ by

$$
\begin{aligned}
\nu_n(C) &= \int_C (2\pi n)^{\frac{d}{2}} \exp\left(-\frac{n\,|\Delta U(t_1,x_1;t_2,x_2)|^2}{2}\right) dt_1 dx_1 dt_2 dx_2 \\
&= \int_C \int_{\mathbb{R}^d} \exp\left(-\frac{|\xi|^2}{2n} + i\langle \xi, \Delta U(t_1,x_1;t_2,x_2)\rangle\right) d\xi\, dt_1 dx_1 dt_2 dx_2.
\end{aligned}
\tag{81}
$$

It follows from Kahane [7] or Testard [12] that (79) will follow if there are positive constants $c_{4,4}$ and $c_{4,5} > 0$ such that

$$\mathbb{E}\big(\|\nu_n\|\big) \geq c_{4,4}, \qquad \mathbb{E}\big(\|\nu_n\|^2\big) \leq c_{4,5}, \tag{82}$$

$$\mathbb{E}\big(\|\nu_n\|_\gamma\big) < +\infty, \tag{83}$$

where $\|\nu_n\| = \nu_n(J)$.

The verifications of (82) and (83) are similar to those in the proof of Theorem 3.3. By Fubini's theorem we have

$$
\begin{aligned}
&\mathbb{E}(\|\nu_n\|) \\
&= \int_J \int_{\mathbb{R}^d} \exp\left(-\frac{|\xi|^2}{2n}\right) \mathbb{E}\exp\left(i\langle \xi, \Delta U(t_1,x_1;t_2,x_2)\rangle\right) d\xi\, dt_1 dx_1 dt_2 dx_2 \\
&= \int_J \int_{\mathbb{R}^d} \exp\left(-\frac{1}{2}\xi\big(n^{-1}\mathrm{I}_d + \mathrm{Cov}(\Delta U(t_1,x_1;t_2,x_2))\big)\xi'\right) d\xi\, dt_1 dx_1 dt_2 dx_2 \\
&= \int_J \frac{(2\pi)^{\frac{d}{2}}}{\sqrt{\det\big(n^{-1}\mathrm{I}_d + \mathrm{Cov}(\Delta U(t_1,x_1;t_2,x_2))\big)}}\, dt_1 dx_1 dt_2 dx_2 \\
&\geq \int_J \frac{(2\pi)^{\frac{d}{2}}}{\sqrt{\det\big(\mathrm{I}_d + \mathrm{Cov}(\Delta U(t_1,x_1;t_2,x_2))\big)}}\, dt_1 dx_1 dt_2 dx_2 := c_{4,4}.
\end{aligned}
\tag{84}
$$

Denote by $\mathrm{Cov}\big(\Delta U(s_1,y_1;s_2,y_2), \Delta U(t_1,x_1;t_2,x_2)\big)$ the covariance matrix of the Gaussian vector $\big(\Delta U(s_1,y_1;s_2,y_2), \Delta U(t_1,x_1;t_2,x_2)\big)$ and let

$$\Gamma = n^{-1}\mathrm{I}_{2d} + \mathrm{Cov}\big(\Delta U(s_1,y_1;s_2,y_2), \Delta U(t_1,x_1;t_2,x_2)\big).$$

Then by the definition of $J$ and (15) in Lemma 1.4, we have

$$
\begin{aligned}
&\mathbb{E}\big(\|\nu_n\|^2\big) \\
&= \int_J\int_J \int_{\mathbb{R}^d}\int_{\mathbb{R}^d} \exp\left(-\frac{1}{2}(\xi,\eta)\,\Gamma\,(\xi,\eta)'\right) d\xi d\eta\, ds_1 dy_1 ds_2 dy_2 dt_1 dx_1 dt_2 dx_2 \\
&= \int_J\int_J \frac{(2\pi)^d}{\sqrt{\det\Gamma}}\, ds_1 dy_1 ds_2 dy_2 dt_1 dx_1 dt_2 dx_2 \\
&\leq \int_J\int_J \frac{(2\pi)^d ds_1 dy_1 ds_2 dy_2 dt_1 dx_1 dt_2 dx_2}{\big[\det\mathrm{Cov}(\Delta U^1(s_1,y_1;s_2,y_2), \Delta U^1(t_1,x_1;t_2,x_2))\big]^{d/2}} \\
&\leq c_{4,6} \int_J\int_J \frac{ds_1 dy_1 ds_2 dy_2 dt_1 dx_1 dt_2 dx_2}{\big[|x_1-y_1| + |x_2-y_2| + |t_1-s_1|^{1/2} + |t_2-s_2|^{1/2}\big]^{d/2}} \\
&\leq c_{4,7} \int_J\int_J \frac{dx_1 dy_1 dx_2 dy_2 dt_1 ds_1 dt_2 ds_2}{\big[|x_1-y_1||x_2-y_2||t_1-s_1||t_2-s_2|\big]^{d/12}} := c_{4,5} < \infty,
\end{aligned}
\tag{85}
$$

where the last inequality follows from $d < 12$. In the above, we have also applied the inequality (36) with $\beta_1 = \beta_2 = 1/6$ and $\beta_3 = \beta_4 = 1/3$.

Similar to (85) and by the same method as in proving (43), we have that $\mathbb{E}(\|\nu_n\|_\gamma)$ is, up to a constant factor, bounded by

$$\int_J \int_J \frac{ds_1 dy_1 ds_2 dy_2 dt_1 dx_1 dt_2 dx_2}{(|x_1-y_1|+|x_2-y_2|+|t_1-s_1|+|t_2-s_2|)^\gamma}$$
$$\times \int_{\mathbb{R}^d}\int_{\mathbb{R}^d} \exp\left(-\frac{1}{2}(\xi,\eta)\,\Gamma\,(\xi,\eta)'\right) d\xi d\eta$$

$$(86) \quad \leq c_{4,8} \int_J \int_J \frac{1}{(|x_1-y_1|+|x_2-y_2|+|t_1-s_1|+|t_2-s_2|)^\gamma}$$
$$\times \frac{dx_1 dy_1 dx_2 dy_2 dt_1 ds_1 dt_2 ds_2}{\left(|x_1-y_1|+|x_2-y_2|+|t_1-s_1|^{1/2}+|t_2-s_2|^{1/2}\right)^{d/2}}$$

$$\leq c_{4,9} \int_0^1 dx_2 \int_0^1 dx_1 \int_0^1 dt_2 \int_0^1 \frac{dt_1}{(t_1^{1/2}+t_2^{1/2}+x_1+x_2)^{d/2}(t_1+t_2+x_1+x_2)^\gamma}$$
$$< \infty,$$

where the last inequality follows from Lemma 2.2, $d < 8$ and the definition of $\gamma$ [We need to consider three cases: $d < 4$, $d = 4$ and $4 < d < 8$, respectively]. This proves (83) and hence Theorem 4.1. □

For $a > 0$ and $b_1 < b_2$, let $K = [a, a+h] \times [b_1, b_1+h] \times [b_2, b_2+h] \subset (0,\infty) \times \mathbb{R}^2$. We choose $h > 0$ small enough, say,

$$h < \frac{b_2 - b_1}{3} \equiv \kappa.$$

Then $|x_2 - x_1| > \kappa$ for all $x_2 \in [b_2, b_2+h]$ and $x_1 \in [b_1, b_1+h]$. We denote the collection of all the cubes $K$ having the above properties by $\mathcal{K}$.

By using Lemma 1.5 and a similar argument as in the proof of Theorem 4.1, we can prove the following dimension result on $L_{\mathrm{II},2}$. We leave the proof to the interested readers.

**Theorem 4.2.** *Let $u = \{u_t(x) : t \geq 0, x \in \mathbb{R}\}$ be a random string process in $\mathbb{R}^d$. If $d \geq 8$, then $L_{\mathrm{II},2} = \emptyset$ a.s. If $d < 8$, then for every $K \in \mathcal{K}$, with positive probability,*

$$(87) \quad \dim_{\mathrm{H}}\left(L_{\mathrm{II},2} \cap K\right) = \dim_{\mathrm{P}}\left(L_{\mathrm{II},2} \cap K\right) = \begin{cases} 3 - \frac{1}{4}d & \text{if } 1 \leq d < 4, \\ 4 - \frac{1}{2}d & \text{if } 4 \leq d < 8. \end{cases}$$

**Remark 4.3.** Rosen [10] studied $k$-multiple points of the Brownian sheet and multiparameter fractional Brownian motion by using their self-intersection local times. It would be interesting to establish similar results for the random string processes.

%beginquote %beginsmall